\documentclass[preprintnumbers,amsmath,amssymb]{revtex4}
\usepackage{color,hyperref}
\usepackage{graphicx}
\usepackage{epsfig}
\usepackage{epstopdf}
\usepackage{ragged2e}
\usepackage{mathrsfs}
\usepackage{csquotes}
\usepackage{tasks}
\usepackage{paralist}
\usepackage{amsthm}
\usepackage[inline]{enumitem}
\usepackage{bm}

\newtheorem{remark}{Remark}[section]
\newtheorem{definition}{Definition}[section]
\newtheorem{theorem}{Theorem}[section]
\newtheorem{proposition}{Proposition}[section]
\newtheorem{lemma}{Lemma}[section]

\begin{document}

\title{Periodic travelling wave solutions of discrete nonlinear Klein-Gordon lattices}

\author{Dirk Hennig}
\author{Nikos I. Karachalios\footnote{Corresponding author. E-mail: karan@uth.gr}}
\affiliation{Department of Mathematics, University of Thessaly, Lamia GR 35100, Greece}

\date{\today}

\begin{abstract}
We prove the existence  of periodic travelling wave solutions for 
general discrete nonlinear Klein-Gordon  systems, considering both cases of hard and soft on-site potentials.
In the case of hard on-site potentials we implement a fixed point theory approach, combining  Schauder's fixed point theorem and the contraction mapping principle. This approach enables us to identify a ring in the energy space for non-trivial solutions to exist, energy (norm) thresholds for their existence and upper bounds on their velocity.  In the case of soft on-site potentials,  the proof of existence of periodic travelling wave solutions  is facilitated by a variational approach based on the  Mountain Pass Theorem. Thresholds on the averaged kinetic energy for these solutions to exist are also derived.
 \end{abstract}

\maketitle
\vspace{0.5cm} \noindent AMS Subject Classifications: 37K40; 37K60; 34C15; 34A33 

\section{Introduction}
Travelling wave solutions (TWSs) in lattice dynamical systems have attracted considerable interest due to their fundamental importance in numerous physical contexts such as the energy transfer in biomolecules, the mobility of dislocations in crystalline  materials and the propagation of pulses in optical systems, to name a few
\cite{Braun2004}-\cite{Satn5}. A variety of fundamental models for the description of such phenomena has been studied, including 
Fermi-Pasta-Ulam-Tsingou, discrete nonlinear Klein-Gordon (DKG) and  discrete nonlinear Schr\"odinger (DNLS) models \cite{Braun2004}-\cite{Kim}, as well as, spatially discrete reaction diffusion systems \cite{Zinner}-\cite{Rowlands} which are relevant to pattern formation resulting from phase transitions. Discrete systems are appropriate  when the scale of decompositions are too small to be effectively described by continuous approximations. 

The existence of solitary travelling wave solutions has been rigorously treated with various approaches from nonlinear analysis \cite{Ambrosetti}-\cite{Zeidler2}, based on direct and minimax variational methods \cite{FWattis}--\cite{PTW1}, reduction to finite dimensional manifolds and normal form techniques \cite{Ioos},\cite{James} and fixed-point methods\cite{Rosen}-\cite{NJF}. For a presentation of several functional-analytic methods and corresponding background for their implementation in nonlinear lattices we refer to the monograph \cite{Pankov1}.

On the other hand, the problem of the existence of {\em periodic TWSs} is much less explored, in particular for second order in time lattice dynamical systems, such as DKG lattices. In the present work, by combining  a  fixed-point approach and variational methods, we study the existence of periodic  TWSs for general DKG systems with anharmonic on-site potentials in  one-dimensional lattices. 

In particular, we study general DKG systems 
described by the following set of coupled oscillator equations
\begin{equation}
\frac{d^2 q_n}{dt^2}=\kappa(q_{n+1}-2q_n+q_{n-1})-V^{\prime}(q_n).
%\,\,\,1 \le n \le N.
\label{eq:system}
\end{equation}
The prime $^{\prime}$ stands for the derivative with respect to $q_n$, the latter being 
the coordinate of the oscillator at site $n$ evolving in an anharmonic on-site potential 
$V(q_n)$.
Each oscillator interacts with its neighbours to the left and right and the 
strength of the interaction 
is regulated by the value of the parameter $\kappa$. 
This system has a Hamiltonian structure related to the energy
\begin{equation}
H=\sum_{n}\left(\frac{1}{2}p_n^2+V(q_n)+
\frac{\kappa}{2}(q_{n+1}-q_{n})^2\right),\label{eq:Ham}
\end{equation}
and it is time-reversible with respect to the involution $p\mapsto -p$.
In the case of a finite lattice, the system (\ref{eq:system}) describes the dynamics of an arbitrary number of $N+1$ oscillators, which are placed equidistantly on the interval $\mathcal{I}=(-L,L)$ of length $2L$. The quantity $\kappa=1/2h^2$  serves as the discretisation parameter, with $h=2L/N$ defining the lattice spacing. The position of the oscillators is given by the discrete spatial coordinate \begin{equation}\label{spatial_cor}
x_n=-L+nh, \quad n= 0,1,2,\ldots,N.
\end{equation}
For finite lattices~\eqref{eq:system} we impose {\it periodic boundary conditions} 
\begin{equation}\label{eq3a} 
q_n = q_{n+N}.
\end{equation}  
We consider TWSs of the form:
\begin{equation}
\label{TWSs}
q_n(t)=Q\left(n-c\,t\right)=Q(z),
\end{equation}
with a $2L-$periodic function $Q(z)$, $z=n-c\,t$, 
where 
$c \in {\mathbb{R}}\setminus \{0\}$ is the velocity.
%are the wave parameters. 
The solutions \eqref{TWSs} satisfy the advance-delay equation
\begin{equation}
c^2{Q}^{\prime\prime}(z)=\kappa (Q(z+1)-2Q(z)+Q(z-1))-V^{\prime}(Q(z)).\label{eq:TW}
\end{equation}
As the system (\ref{eq:system}) possesses the time-reversibility symmetry it suffices to consider  $c>0$.

The presentation of the results is as follows: 
Section \ref{section:hard} deals with the existence of periodic TWSs in finite lattices with {\em hard on-site potentials}.
Imposing periodic boundary conditions, the existence of periodic TWSs is proved 
%by an up to our knowledge, novel 
by a fixed-point approach, based on  Schauder's fixed point  theorem and contraction mapping principles. The Schauder's fixed point method was made use of in our recent works \cite{JFPT}, \cite{DNBRS}  in order to establish the existence of breather solutions. To our knowledge, the extension of this approach to the problem of the existence of TWs is novel, and as in the case of breathers, can offer  a comparatively concise method (if compared to the aforementioned ones), for treating the advance-delay equation \eqref{eq:TW} in the suitable functional set-up of Sobolev spaces of periodic functions. The combination with the contraction mapping argument yields the identification of a ring in the energy space where nontrivial solutions  exist, together with the derivation of energy thresholds for their existence, and an upper bound for their velocity. For the latter, the contraction mapping argument is applied on a map defined by a suitable auxiliary linear, non-homogeneous problem stemming from \eqref{eq:TW}, whose solvability is guaranteed by the Friedrichs extension theorem. We should remark that while the current literature on the energy threshold problem concerns mainly breather solutions for DNLS models \cite{ETH1},\cite{ETH2} (see also \cite{Ch1},\cite{Ch2},\cite{NJF},\cite{ETH3} and references therein), to the best of our knowledge,  the results of the present work are the first on energy thresholds for TWSs. 

In section \ref{section:soft} we consider the case of 
{\em soft on-site potentials}. This time,  for the proof of existence of non-trivial periodic TWSs, we implement a variational approach based on the Mountain Pass Theorem. Still in this case, threshold values for the average kinetic energy of TWSs are derived with the aid of the contraction mapping principle. 

In both cases of potentials, the estimates on the norms of the solutions and the various thresholds are derived under certain conditions exhibiting a coherent dependence  on the lattice parameters, the frequency of the TWs and their velocity which is proved to be physically relevant (e.g. see the concluding remarks of section II for the case of the hard potentials). 
 
\section{Hard on-site potentials}\label{section:hard}
\setcounter{equation}{0}
In this section we consider  hard on-site (uniformly convex) potentials $V(x)$  assuming that the following assumption holds on them:

\vspace*{0.5cm}

\noindent {\bf{A:}}  
$V:{\mathbb{R}}\rightarrow {\mathbb{R}}$ is  
non-negative and at least twice continuously 
differentiable and is characterised by the following properties:
\begin{eqnarray}
 V(0)&=&V^\prime(0)=0,\,\,\,V^{\prime \prime}(0)>0,\label{eq:potential0}\\
V(x)&=&V(-x),\;\;\mbox{for all\  $x\in\mathbb{R}$},\label{eq:evpotential}
\end{eqnarray}
The unique equilibrium at $x=0$ is a global minimum of $V(x)$. 
Further, we assume that $V$ satisfies for some positive constants $\overline{m}$ $\alpha$, $\beta$, and $K$,  the conditions 
\begin{eqnarray}
|V''(x)|&\leq& \overline{m}|x|^\alpha,\;\;\forall x\in\mathbb{R},
	\label{eq:assumptions1}\\
	|V^\prime(x_1)-V^\prime(x_2)|&\le& K(|x_1|^\beta+|x_2|^\beta)|x_1-x_2|.
	\label{eq:assumptions2}	
\end{eqnarray}

We consider periodic solutions $Q(z)=Q(z+2L)$ of zero-mean satisfying
\begin{equation}
\int_{-L}^LQ(z)dz=\int_{-L}^L V'(Q(z))dz=0,\label{eq:zeromean}
\end{equation}
which is satisfied for even potentials \eqref{eq:evpotential} and periodic functions exhibiting the symmetry 
\begin{equation}
	\label{eq:symm}
 Q(z-L/2)=-Q(L/2-z).
\end{equation}
We shall further discuss the conditions \eqref{eq:evpotential}, \eqref{eq:zeromean} and \eqref{eq:symm} of the assumption \textbf{A} below (see Theorem \ref{ThII.1}). 
Under assumption \textbf{A},  we reformulate the original problem  (\ref{eq:TW}) as a fixed point
problem for a suitably defined operator in a Banach space. Motivated by our recent works \cite{JFPT},\cite{DNBRS} applying Schauder's Fixed Point Theorem \cite{Zeidler} to prove  the existence of breathers in discrete nonlinear KG lattices,  we show that the same Theorem  is applicable, in order to prove the existence of periodic travelling waves. 

In particular, we will use the following version of Schauder's Fixed Point Theorem (see e.g. in \cite{Zeidler}): 

\begin{theorem}
 \label{theorem:SFPT}
 Let $G$ be a non-empty closed convex subset of the 
Banach space $X$. Suppose that
$U:\,\,G\rightarrow\,G$  is a continuous and compact  map.  
then  $U$ has a fixed point in $G$.
\end{theorem}
%%%%%
\subsection{Existence of periodic TWSs}
We start by introducing appropriate function spaces on which our methods will be employed. We consider first the space
\begin{equation}
{{\mathcal{X}}}_0\,=\,\left\{\,Q\in L^{2}_{per}(-L,L)\,\vert\,\,\,\int_{-L}^LQ(s)ds=\int_{-L}^L V'(Q(s))ds=0,\,\,\,Q(s-L/2)=-Q(L/2-s)\right\}\,,
\end{equation}
that is, the space of $L^2_{per}(-L,L)$, $2L-$periodic, square integrable functions of zero mean, on the account of conditions \eqref{eq:evpotential}, \eqref{eq:zeromean} and \eqref{eq:symm}.  The space $\mathcal{X}_0$ is endowed with the norm,
\begin{equation}
|| Q||_{{{\mathcal{X}}}_0}=\left(
\frac{1}{2L}\,
\int_{-L}^L\,Q(s)^2ds\right)^{1/2}= \left(\sum_{k\in {\mathbb{Z}}\setminus \{0\}}|\hat{{Q}}_{k}|^2 dk\right)^{1/2}.
\end{equation}
$\hat{Q}_{k}$ determines the Fourier-coefficient in the Fourier series expansion:
\begin{eqnarray}
Q(s)&=&\sum_{k\in {\mathbb{Z}}\setminus \{ 0 \}}
\hat{Q}_{k} \exp\left( i\Omega k s\right),\,s\in (-L,L),\,\,\,\Omega=\frac{\pi}{L},\label{eq:Fourier}\\
{\hat{Q}}_{k}&=&\frac{1}{2L} \int_{-L}^{L} Q(s)\exp(-i \Omega k s)ds,
\,\,\,\hat{Q}_{k}=\overline{\hat{Q}}_{-k},\label{eq:Fourier1}
\end{eqnarray}
of $Q(s)$ and $\overline{x}$ denotes the complex conjugate of 
$x$.
%%%%%%
Evidently, $\mathcal{X}_0$ is a closed subspace of $L^{2}_{per}(-L,L)$. 
We shall also use the spaces
\begin{equation}
\label{sp1}
{{\mathcal{X}}}_1\,=\,\left\{\,Q\in H^{1}_{per}(-L,L)\,\vert\,\,\,\int_{-L}^LQ(s)ds=\int_{-L}^L V'(Q(s))ds=0,\,\,\,Q(s-L/2)=-Q(L/2-s)\,\right\}\,,
\end{equation}
endowed with the scalar product and induced norm
\begin{eqnarray}
	\label{homX1}
	(U,Q)_{\mathcal{X}_1}=\int_{-L}^{L}U'(s)Q'(s)d{s},\;\;||Q||^2_{\mathcal{X}_1}=\int_{-L}^{L}Q'(s)^2d {s}=\sum_{k\in {\mathbb{Z}}\setminus \{0\}}(k\Omega)^2
		|\hat{{Q}}_{k}|^2.
\end{eqnarray} 
It can be easily checked that the Poincar\'{e} inequality
\begin{eqnarray}
\label{e6}
\int_{-L}^{L}Q(s)^2ds\leq {C}(L)\int_{-L}^{L}Q'(s)^2ds,
\end{eqnarray}
holds with  the optimal constant ${C}(L)={L^2}/{\pi^2}$ when the Hilbert space $\mathcal{X}_1$ is used. Note that the validity of the  Poincar\'{e} inequality \eqref{e6} allows for the usage of the homogeneous norm \eqref{homX1}. 
Of use will also be the space
\begin{equation}
\label{sp2}
{{\mathcal{X}}}_2\,=\,\left\{\,Q\in H^{2}_{per}(-L,L)\,\vert\,\,\,\int_{-L}^LQ(s)ds=\int_{-L}^L V'(Q(s))ds=0,\,\,\,Q(s-L/2)=-Q(L/2-s)\,\right\}.
\end{equation}
Due to the Poincar\'{e} inequality \eqref{e6}, for ${{\mathcal{X}}}_2$ we shall use the homogeneous norm
\begin{eqnarray}
|| Q ||_{\mathcal{X}_2}^2=\int_{-L}^{L}Q''(s)^2d {s}=
\sum_{k\in {\mathbb{Z}}\setminus \{0\}}(k\Omega)^4
|\hat{{Q}}_{k}|^2.\label{eq:norml2}
\end{eqnarray} 
It can be readily seen that $\mathcal{X}_{2}$  is a closed subspace of $H_{per}^{2}(-L,L)$.
%and we used the notation $H^0=L^2$. 
Furthermore,  $\mathcal{X}_2$ is compactly embedded in $\mathcal{X}_0$ ($\mathcal{X}_2\Subset \mathcal{X}_{0}$).  

In some cases, we will facilitate variational methods. In particular,
solutions of (\ref{eq:TW}) will be considered as critical points of the action functional  $S:\,\mathcal{X}_1 \rightarrow {\mathbb{R}}$ given by
\begin{equation}
\label{AcFun}
S(Q)=\int_{-L}^L\left[\frac{1}{2}\left(cQ^\prime(z)\right)^2-V(Q(z))-
\frac{1}{2}\kappa\,[Q(z+1)-Q(z)]^2\right]dz.\nonumber
%\label{eq:action}
\end{equation}
Related with the action functional is the (total) energy functional $E$ on $\mathcal{X}_1$:
\begin{equation}
E(Q)=
%\int_{-L}^L\, e(cQ(z),Q^\prime(z))du=
\int_{-L}^L\left[\frac{1}{2}\left(cQ^\prime(z)\right)^2+V(Q(z))+
\frac{\kappa}{2}\,[Q(z+1)-Q(z)]^2\right]dz,\label{eq:Ham1}
\end{equation}
which will be also involved in the derivation of several bounds for the solutions.

Let us denote by $C_{0,*}$ the constant of the embedding $\mathcal{X}_2\Subset\mathcal{X}_0$, and  by $C_{2,*}$ the constant of the embedding $\mathcal{X}_2\subset L^{\infty}([-L,L])$.
We now proceed  to the statement and proof of the result on the existence of periodic TWSs:
\begin{theorem}
	\label{ThII.1}	
	Let  assumption {\bf{A}} hold.  
	If
	\begin{equation}
		% \frac{1+\Omega^2}{\,\Omega^2-4\sum_{j=1}^{N_c}\,
			%  \kappa_j\sin^2\left(\frac{\Omega }{2}jk\right)}{\overline{m}}<1
		%\overline{m}\le \frac{1-{\cal S}}{1+\Omega^{-2}},\,\,\,{\cal S}<1,\label{eq:As}
		%\overline{m}\le \Omega^2(1-{\cal S}),
		\frac{4\kappa}{c^2}<\Omega^2,
		%\,\,\,{\cal S}<1,
		\label{eq:As}
	\end{equation}
	%  where 
	%  \begin{equation}
		%   {\cal S}=\sup_{x\in {\mathbb{R_+}}} g(x)<1,
		%  \end{equation}
	%  with the  
	%  function $g\in C({\mathbb{R_+}},{\mathbb{R}})$ (\,${\mathbb{R}}_+=[0,\infty)$\,) given by
	%  \begin{equation}
		%  g(x)=\frac{2}{x}\kappa\,\sin^2(x), 
		% \end{equation}
	then 
	there exists a periodic TWS 
	$q_n(t)=Q\left(n-\,ct\right)\equiv Q(z)$ of (\ref{eq:system}) with 
	$Q \in H^2(-L,L)$ and 
	\begin{equation}
		\label{B1}
		||Q||_{\mathcal{X}_0}\le \left( \frac{c^2\Omega^2-4
			\kappa}{{K}C_{3,*}^{\beta}}\right)^{1/\beta}:=\mathscr{R}_{max},
	\end{equation}
where the constant $C_{3,*}$ is given by
\begin{equation}
	\label{newc}
C_{3,*}:=\frac{C_{2,*}}{C_{0,*}},	
\end{equation}
such that 
	%for some $T=2\pi/\Omega>0$ 
	\begin{equation}
		Q(z+2L)=Q(z),\,\,\forall z \in {\mathbb{R}}.
	\end{equation}
\end{theorem}
%\vspace*{0.5cm}

\noindent{\bf Proof:} We seek    
$2L-$periodic functions $Q\in H^{2}(-L,L)$ satisfying (\ref{eq:zeromean}-\eqref{eq:symm} and solve Eq.\,(\ref{eq:TW}). For the following discussions it is suitable to re-write Eq.\,(\ref{eq:TW}) as:
\begin{equation}
	{Q}^{\prime \prime}(z)-\frac{\kappa}{c^2} (Q(z+1)-2Q(z)+Q(z-1))=-\frac{1}{c^2}V^{\prime}(Q(z)).
	\label{eq:ref1}
\end{equation}
Thus only the right-hand side of  (\ref{eq:ref1}) features terms nonlinear in $Q$. Ultimately, 
(\ref{eq:ref1}) shall be expressed as a fixed point equation in $Q$. 

To  proceed, we recall first some basic facts for the intersection of Banach spaces \cite[Lemma 2.3.1 and Theorem 2.7.1]{interB}: For two Banach spaces $X$ and $Y$, the intersection $X\cap Y$ is a Banach space endowed with the norm
\begin{eqnarray}
\label{intera}
||x||_{X\cap Y}=||x||_{X}+||x||_{Y},\;\;\mbox{for all $x\in X\cap Y$.}
\end{eqnarray}
Clearly, from \eqref{intera}, we have that 
\begin{eqnarray}
	\label{interb}
||x||_{X}\leq ||x||_{X\cap Y},\;\;\mbox{and}\;\;||x||_{Y}\leq 	||x||_{X\cap Y},\;\;\mbox{for all $x\in X\cap Y$.}
\end{eqnarray}  
Moreover, in the particular case where $Y$ is continuously embedded in $X$ with an embedding constant $c_0$, that is $||x||_{x}\leq c_0||x||_{Y}$, for all $x\in Y$, then, due to \eqref{intera} and \eqref{interb},
\begin{eqnarray}
	\label{interc}
||x||_{X}\leq ||x||_{X\cap Y}\leq c_0^*||x||_{Y},\;\;\mbox{for all $x\in Y$,}\;\;c_0^*=1+c_0.
\end{eqnarray}
Hence, if $x$ is in a closed ball of $Y$ of radius $r$, then due to \eqref{interc}, $x$ is in the closed ball of $X\cap Y$ of radius $c_0^{*}r$, and in the closed ball of the same radius of $X$.  With these preparations we apply the above for $X=\mathcal{X}_0$ and $Y=\mathcal{X}_2$, and   
we define a  convex set of $\mathcal{X}_0\cap\mathcal{X}_2$, as follows: First, we consider arbitrary elements $P\in\mathcal{X}_0\cap\mathcal{X}_2$, such that 
\begin{equation}
	\label{cset1}
||P||_{\mathcal{X}_2}\leq \varrho.	
\end{equation}	
Since $P\in\mathcal{X}_2\cap\mathcal{X}_0$ and $\mathcal{X}_2\Subset\mathcal{X}_0$  with compact embedding, it holds by application of \eqref{interc}, that
\begin{eqnarray}
	\label{cset2}
	||P||_{\mathcal{X}_0}\leq ||P||_{\mathcal{X}_0\cap\mathcal{X}_2}\leq C_{0,*}||P||_{\mathcal{X}_2}\leq C_{0,*}\varrho,
\end{eqnarray}
Then, we set for simplicity,
\begin{eqnarray}
	\label{cset3}
	R=C_{0,*}\varrho,
\end{eqnarray}
and we consider the convex set of $\mathcal{X}_0\cap\mathcal{X}_2$ as
\begin{equation}
	{\mathcal{Y}}_0\,=\,\left\{\,P\in \mathcal{X}_0\cap\mathcal{X}_2\,\,:
	\,\,|| P||_{\mathcal{X}_0}\le R\,\right\}.
\end{equation}

Evidentially, the set $\mathcal{Y}_0$ is well defined and non-empty due to \eqref{cset1} and \eqref{cset2}. As we will prove below, $R$ (and thus, $\varrho$, due to \eqref{cset3}), will be explicitly determined, giving rise to the estimate \eqref{B1} for the TWs. 

We relate the left-hand side  of (\ref{eq:ref1}) to  the 
linear mapping: $M\,:\,\mathcal{X}_2\,\rightarrow\,\mathcal{X}_{0}$:
\begin{equation}
	M(P)={P}^{\prime \prime}(z)-\frac{\kappa}{c^2}(P(z+1)-2P(z)+P(z-1)).  
\end{equation}
As a next step, we establish the invertibility of this mapping and derive an upper  bound 
for the norm of its inverse. 

Applying the operator $M$ to the Fourier elements $\exp(i\,\Omega ks)$ in (\ref{eq:Fourier}), results in
\begin{equation}
	M\exp(i\,\Omega ks)=\nu_k \exp(i\,\Omega ks),
\end{equation}
where 
\begin{equation}
	\nu_k= -
	\Omega^2\,k^2+4\,
	\frac{\kappa}{c^2}\sin^2\left(\frac{\Omega }{2}k\right).
\end{equation}
The hypothesis (\ref{eq:As}) guarantees that $\nu_k \ne 0$, for all $k \in {\mathbb{Z}}\setminus\{0\}$. 
Therefore, the mapping $M$ possesses an inverse obeying 
$M^{-1}\exp(i\,\Omega ks)=(1/\nu_k)\exp(i\,\Omega ks)$. 
%%%%%%%%%%%%%Newton
For the norm of the 
linear operator $M^{-1}:\,\,\mathcal{X}_{0} \rightarrow \mathcal{X}_2$, one gets, by using the hypothesis  (\ref{eq:As}), the upper bound:
\begin{eqnarray}
	\lefteqn{ || M^{-1} ||_{\mathcal{X}_{0},\mathcal{X}_2}
		%\equiv || M^{-1}||
		=\sup_{0 \neq Q \in \mathcal{X}_0}\frac{|| M^{-1}\,Q ||_{\mathcal{X}_2}}
 		{|| Q ||_{\mathcal{X}_0}}}\nonumber\\
	&=&\sup_{0 \neq Q \in \mathcal{X}_0}
	\frac{\left( \sum_{l}^\prime (\Omega\,l)^4\left|\frac{1}{\nu_l} 
		\hat{Q}_{l}\exp(i\,l\Omega s)\right|^2 \right)^{1/2}}{|| Q ||_{\mathcal{X}_0}}\nonumber\\
	%%%%%%%%%%%%%%%%%%% 
	&\le & \sup_{l\in {\mathbb{Z}}\setminus \{ 0 \}} \frac{\sqrt{(\Omega\,l)^4} }{|\nu_l|}\,
	\cdot \sup_{0 \neq Q \in \mathcal{X}_0}\frac{\left(\sum_{l}^\prime 
		|\hat{Q}_{l}|^2\right)^{1/2}}{|| Q ||_{\mathcal{X}_0}}
	\nonumber\\
	&=&\sup_{l\in {\mathbb{Z}}\setminus \{ 0 \}} \frac{{(\Omega\,l)^2}}{|-
		\Omega^2\,l^2+\frac{4\kappa}{c^2}\sin^2\left(\frac{\Omega }{2}l\right)|}\nonumber\\
	&\le& 
	\frac{\Omega^2}{\,\Omega^2-\frac{4\kappa}{c^2}},\label{eq:boundL}
\end{eqnarray}
which verifies the boundedness of $M^{-1}$;  note that we have used the notation $\sum_l^\prime = \sum_{l\in {\mathbb{Z}}\setminus \{ 0 \}}$. Obviously, the same estimate \eqref{eq:boundL} verifies the boundedness of $M^{-1}$ as an operator $M^{-1}:\mathcal{X}_2\cap\mathcal{X}_0\rightarrow\mathcal{X}_0$ On the other hand, we remark for later use, that if we  consider only the lower-order terms of the estimate \eqref{eq:boundL},  we get that
\begin{eqnarray}
	\label{x0x0}
	||M^{-1}||_{\mathcal{X}_0,\mathcal{X}_0}\leq  \frac{c^2}{c^2\Omega^2-4\kappa}.
\end{eqnarray}

Associated with the right-hand side of (\ref{eq:ref1}) we introduce 
the nonlinear operator $N\,:\,\mathcal{Y}_{0}\,\rightarrow\,\mathcal{Y}_{0}$, as
\begin{equation}
	N(P)= -\frac{1}{c^2}V^{\prime}(P).\label{eq:G}
\end{equation}

From the definition of the operator $N$ in \eqref{eq:G}, the necessity for the conditions \eqref{eq:evpotential},\eqref{eq:zeromean} and \eqref{eq:symm} becomes evident, as they guarantee the closeness of the operator $N$ in $\mathcal{Y}_0$. To verify continuity of the  operator $N$  on $\mathcal{X}_{0}$, for any $P\in \mathcal{Y}_0$, (which implies by the definition of $\mathcal{Y}_0$ that $P\in \mathcal{X}_{0}\cap\mathcal{X}_2$), 
we prove that  $N$ is Frechet differentiable with a locally bounded derivative. Using (\ref{eq:assumptions1})  we have
\begin{equation}
	N^{\prime}(P):\,h\in \mathcal{X}_{0}\,\mapsto N^{\prime}(P)[h] =-\frac{1}{c^2}V^{\prime\prime}(P)h \in \mathcal{X}_{0}.
\end{equation}
Then, since $P\in\mathcal{X}_2\cap\mathcal{X}_0$ and $\mathcal{X}_2\subset L^{\infty}([-L,L])$ with continuous embedding, we have that
\begin{eqnarray}
	|| N^{\prime}(P)[h]||_{\mathcal{X}_{0}}&=&
	|| -\frac{1}{c^2}V^{\prime \prime}(P)\,h ||_{\mathcal{X}_0}\nonumber\\
	&=& \frac{1}{c^2}\left(\frac{1}{2L}\int_{-L}^L| V^{\prime\prime}(P(z))\,h(z)|^2dz\right)^{1/2}\le \frac{1}{c^2}\left(\overline{m}^2 \max_{-L\le z\le L}|P(z)|^{2\alpha}\frac{1}{2L}\int_{-L}^L|h(z)|^2dz\right)^{1/2}\nonumber\\
	&\le& \frac{\overline{m}}{c^2}\,\max_{-L\le z\le L}|P(z)|^{\alpha}\,|| h||_{\mathcal{X}_{0}}\nonumber\\
	&\leq& \frac{\overline{m}}{c^2}C_{2,*}^{\alpha}||P||_{\mathcal{X}_2}^{\alpha}|| h||_{\mathcal{X}_{0}}\leq A(\overline{m},c^2,C_{2,*},\varrho)|| h||_{\mathcal{X}_{0}},
\end{eqnarray}
for the constant
$$A=A(\overline{m},c^2,C_{2,*},\varrho)=\frac{\overline{m}}{c^2}C_{2,*}^{\alpha}\varrho^{\alpha},$$ 
proving the local boundedness of the differential. Then, one concludes that the Frechet derivative is locally bounded as
\begin{equation}
	|| N^{\prime}(P)||_{{\cal{L}}(\mathcal{X}_{0},\mathcal{X}_{0})}\le A.\label{eq:Frechet}
\end{equation} 
With the aid of the  locally bounded derivative, we can prove the local Lipschitz continuity of $N$ as follows:
\begin{eqnarray}
	|| N(P)-N(Q) ||_{\mathcal{X}_{0}}&\le& \sup_{S \in [P,Q]} || N^{\prime}(S)||_{{\cal{L}}(\mathcal{X}_{0},\mathcal{X}_{0})}|| P-Q||_{\mathcal{X}_{0}}\nonumber\\
	&\le&A|| P-Q||_{\mathcal{X}_{0}}.\label{eq:LC}
\end{eqnarray}
As the range of $N$ is concerned, by using condition  \eqref{eq:assumptions2}, we derive that for any $P\in \mathcal{Y}_0$,
\begin{eqnarray}
	|| N(P) ||_{\mathcal{X}_{0}}=\frac{1}{c^2}|| \,-V^{\prime}(P)\, ||_{\mathcal{X}_0}&=&\frac{1}{c^2}\left(\frac{1}{2L}\int_{-L}^L|-V^{\prime}(P(z))|^2du\right)^{1/2}\le \frac{1}{c^2}\left(\frac{1}{2L}\int_{-L}^L K^2| P(z)|^{2(\beta+1)}dz\right)^{1/2}\nonumber\\
	& \le & \frac{K}{c^2}\,\max_{-L\le z \le L}|P(z)|^{\beta}|| P ||_{\mathcal{X}_{0}}\nonumber\\
	& \le &
	\frac{K}{c^2}C_{2,*}^{\beta}||P||_{\mathcal{X}_2}^{\beta}||P||_{\mathcal{X}_0}\nonumber\\
	& \le &
	\frac{K}{c^2}C_{2,*}^{\beta}C_{0,*}||P||_{\mathcal{X}_2}^{\beta+1}\nonumber\\
	&\leq &\frac{K}{c^2}C_{2,*}^{\beta}C_{0,*}\varrho^{\beta+1}
	= \frac{K}{c^2}\left(\frac{C_{2,*}}{C_{0,*}}\right)^{\beta}R^{\beta+1}.
	\label{eq:rangeN}
\end{eqnarray}
Thus, we proved in \eqref{eq:rangeN}, that for any $P\in\mathcal{Y}_0$,
\begin{eqnarray}
	\label{finesr}
|| N(P) ||_{\mathcal{X}_{0}}\leq \frac{KC_{3,*}^{\beta}}{c^2}R^{\beta+1},
\end{eqnarray}
with the constant $C_{3,*}=\frac{C_{2,*}}{C_{0,*}}$, as defined in \eqref{newc}.

At last, we express the problem  (\ref{eq:ref1})
as a fixed point equation in terms of a mapping $\mathcal{Y}_0\,\rightarrow\, \mathcal{Y}_{0}$:
\begin{equation}
	Q=M^{-1}\circ N(Q)\equiv \mathscr{L}(Q).\label{eq:compose}
\end{equation}
Using (\ref{x0x0}) and (\ref{eq:rangeN}), we have
\begin{eqnarray}
	||\mathscr{L}(Q)||_{\mathcal{X}_0}&=&|| M^{-1}(\,N(Q))||_{\mathcal{X}_0}
	\le || M^{-1} ||_{\mathcal{X}_0,\mathcal{X}_0}||  \,N(Q)||_{\mathcal{X}_0}\nonumber\\
	& \le& 
	\frac{{K}C_{3,*}^{\beta}}{c^2\Omega^2-4
		\kappa}R^{\beta+1}\le R,
\end{eqnarray}
assuring by assumptions (\ref{eq:As}) and \eqref{B1}, that indeed,
\begin{equation}
	\mathscr{L}(\mathcal{Y}_0)\subseteq \mathcal{Y}_0.
\end{equation} 
Furthermore, since it holds that $\mathscr{L}(Q)\in \mathcal{X}_2$ for all $Q \in \mathcal{Y}_0\subseteq \mathcal{X}_0$, one has $\mathscr{L}(\mathcal{Y}_0)\subseteq \mathcal{X}_2 \cap \mathcal{Y}_0$, 
and as the embedding of $\mathcal{X}_2$ in $\mathcal{X}_0$ is compact, the operator $\mathscr{L}$ is compact.  
It remains to prove that $\mathscr{L}$ is continuous on $\mathcal{Y}_{0}$: 
For arbitrary $P_1,P_2\in \mathcal{Y}_{0}$, we have
\begin{eqnarray}
	||\mathscr{L}(P_1)-\mathscr{L}(P_2)||_{\mathcal{X}_0}&=&|| M^{-1}(N(P_1))- M^{-1}(N(P_2))
	||_{\mathcal{X}_2}\le 
	|| M^{-1}||_{\mathcal{X}_{0},\mathcal{X}_{0}}
	\,|| N(P_1)- N(P_2)||_{\mathcal{X}_{0}}\nonumber\\
	&\le& \frac{Ac^2}{c^2\Omega^2-4\kappa} ||P_1-P_2||_{\mathcal{X}_0}<\epsilon,
\end{eqnarray}
if 
\begin{equation}
 || P_1-P_2||_{\mathcal{X}_0} < \delta=\frac{c^2\Omega^2-4\kappa}{Ac^2}\epsilon,
\end{equation}
for any given $\epsilon>0$, verifying that $\mathscr{L}(Q)$ is continuous on $\mathcal{Y}_0$. Thus, all the assumptions of Schauder’s fixed point theorem are satisfied, and hence, the fixed point equation $Q =  \mathscr{L}(Q)$ has at least one solution. 
\ \ $\square$

\begin{remark} (Regularity of travelling waves). By the Sobolev embeddings, the obtained $H^2$-travelling wave solutions $Q$ are $C^{1}$. Therefore, due to Eq.\,(\ref{eq:TW}) it holds that 
${Q}^{\prime \prime} \in C^{1}$ and conclusively $Q\in C^{3}$, that is, they  are classical solutions.
\end{remark}
\subsection{Existence of non-trivial  TWSs and an energy threshold}
In this section, we consider the problem of existence and non-existence of non-trivial TWSs with frequencies satisfying the 
strengthened condition
\begin{eqnarray}
\label{eq:As1}
\Omega^2>\frac{4\kappa+KC_{3,*}^{\beta}}{c^2},
\end{eqnarray}
if compared to \eqref{eq:As}. The motivation for assuming \eqref{eq:As1}, is explained in Theorem \ref{ThII.2} concluding the section.  We start by stating  a result on non-existence of non-trivial TWSs with frequencies satisfying \eqref{eq:As1}.
\begin{proposition}
\label{PropII.1}	
Suppose that conditions \textbf{A}  and \eqref{eq:As1} hold and that
 \begin{equation}
||Q||_{\mathcal{X}_0}\le R< 
\left(\frac{c^2\Omega^2-4
 \kappa}{{K}C_{3,*}^{\beta}}\right)^{1/(1+\beta)}:=\mathscr{R}_{crit}
%\frac{c^2\Omega-4\kappa}{\overline{m}C_{2,*}^{\beta}}
.\label{eq:uni}
 \end{equation}
 Then the equation $Q=\mathscr{L}(Q)$ has only the trivial solution, that is,  there exist no non-trivial TWSs. 
\end{proposition} 
 {\bf Proof:} If (\ref{eq:uni}) holds, we get
 $|| \mathscr{L}(Q) ||_{\mathcal{X}_0}=
 || M^{-1}(N(Q))||_{\mathcal{X}_0} \le 
 || M^{-1}||_{\mathcal{X}_0,\mathcal{X}_0}\cdot || N(Q)||_{\mathcal{X}_0}<1$.
 Thus $\mathscr{L}$ is a  contraction, 
 and the Contraction Mapping Theorem implies that there is a
 unique function $Q$ that solves the equation $Q=\mathscr{L}(Q)$. Since $\mathscr{L}(0)=0$, this unique solutions is the trivial one. 

 Note that due to (\ref{eq:potential0}) the potential $V(x)$ can be expressed as $V(x)=(1/2)x^2+W(x)$,  
  with
 $W(x)\ge 0$ for all $x\in \mathbb{R}$, satisfying \eqref{eq:assumptions1} and \eqref{eq:assumptions2} (see also \cite[Eq. (2.7) and Eq. (2.8)]{DNBRS}).
Furthermore, let us observe that
as the (conserved) energy functional (\ref{eq:Ham1}) is coercive, it can be used to bound the $L^2(-L,L)$-norm of the TWSs as follows 
 \begin{eqnarray}
 E(Q)&=&
 \int_{-L}^L\left[\frac{1}{2}\left(Q^\prime(z)\right)^2+\frac{1}{2}Q^2(z)+W(Q(z))+
 \frac{\kappa}{2}\,[Q(z+1)-Q(z)]^2\right]dz\nonumber\\
 &=& \int_{-L}^L\left[\frac{1}{2}\left(Q^\prime(z)\right)^2+W(Q(z))+
 \frac{\kappa}{2}\,[Q(z+1)-Q(z)]^2\right]dz+\frac{1}{2}||Q||_{L^2(-L,L)}^2,
 \end{eqnarray}
 so that 
 \begin{equation}
 || Q||_{\mathcal{X}_0}\le \sqrt{2E}.\label{eq:estimatesH1}
 \end{equation}
 In conclusion, if the energy of the lattice system is less than $E_{crit}(c,\kappa,K,\Omega)=\sqrt{2R_{crit}}$, no TWSs of given values for $c,\kappa, K,\Omega$ exist. 
 \ \ $\square$

\begin{figure}[tbp!]
	\begin{center}
		\begin{tabular}{cc}
			\includegraphics[width=.42\textwidth]{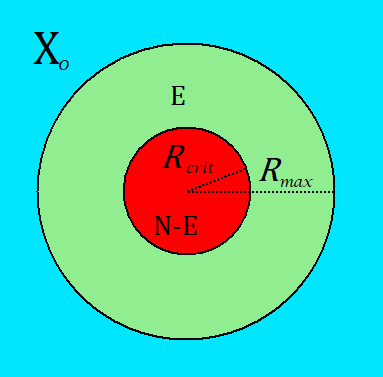}		
		\end{tabular}
	\end{center}
	\caption{Illustration  of the statement of Theorem \ref{ThII.2}:  Non-trivial travelling wave solutions exist in the ring  $\mathbf{E}=\overline{\mathcal{B}(0,\mathscr{R}_{max})}\setminus \overline{\mathcal{B}(0,\mathscr{R}_{crit})}$ of the space $\mathcal{X}_0$ [light (green) coloured area ($\mathbf{E}$)]. The radius $\mathscr{R}_{crit}$ of the closed ball $\overline{\mathcal{B}(0,\mathscr{R}_{crit})}$ follows from the contraction mapping Theorem (see Proposition \ref{PropII.1}) and defines a threshold for non-existence of non-trivial solutions [in the darker (red) coloured area ($\mathbf{N}-\mathbf{E})$]. The radius $\mathscr{R}_{max}$ of the closed ball $\overline{\mathcal{B}(0,\mathscr{R}_{max})}$ results from  the Schauder's fixed point theorem (see Theorem \ref{theorem:SFPT}).}
	\label{fig1}
\end{figure}

Combining the existence Theorem \ref{ThII.1} and Proposition \ref{PropII.1} we identify a ring in the phase space $\mathcal{X}_0$ with quantified radii in which the non-trivial travelling wave solutions with frequencies satisfying the enhanced condition \eqref{eq:As1}, exist. This is illustrated in the cartoon of Fig. \ref{fig1}.
\begin{theorem}
\label{ThII.2}
Let the assumption {\bf A} and the condition \eqref{eq:As1} hold.
The system (\ref{eq:system})  possesses non-trivial TWSs only if 
\begin{equation}
\label{ring1}
\mathscr{R}_{crit}= \left(\frac{c^2\Omega^2-4
	\kappa}{{K}C_{3,*}^{\beta}}\right)^{1/(1+\beta)}\leq ||Q||_{\mathcal{X}_0}\leq 
\mathscr{R}_{max}=\left(\frac{c^2\Omega^2-4
	\kappa}{{K}C_{3,*}^{\beta}}\right)^{1/\beta}.
\end{equation}
\end{theorem}
\textbf{Proof:} On the one hand, Theorem \ref{ThII.1} establishes the existence of solutions when $||Q||_{\mathcal{X}_0}\leq \mathscr{R}_{max}$. On the other hand, according to Proposition \ref{PropII.1}, if $||Q||_{\mathcal{X}_0}< \mathscr{R}_{crit}$ only the trivial solution exists. Thus, a non-trivial solution exists only if 
\begin{eqnarray*}
\mathscr{R}_{crit}<||Q||_{\mathcal{X}_0}\leq 	\mathscr{R}_{max},
\end{eqnarray*}
that is, when \eqref{ring1} is satisfied. 
For the latter to be valid, we require  $\mathscr{R}_{crit}<\mathscr{R}_{max}$, motivating the condition \eqref{eq:As1} on the frequencies of the TWSs. \ \ $\square$
\subsection{Upper bound for the velocity}\label{subsection:upperbound}
\label{SECIIC}
By using a fixed point approach, we derive upper bounds for the velocity $c$ of TWSs. 
For this purpose, we rewrite equation \eqref{eq:TW} in the following operator form:
\begin{eqnarray}
\label{eop3}
-Q''(z)={\frac{1}{c^2}\left\{\kappa\left([Q(z)-Q(z-1)]-[Q(z+1)-Q(z)]\right)+V^\prime(Q(z))\right\}}.
\end{eqnarray}
We recall some basic auxiliary results, starting with the Friedrichs extension Theorem.
\begin{theorem}
	\label{T1}
	Let $\mathcal{L}_0:D(\mathcal{L}_0)\subseteq {X}_0\rightarrow {X}_0$ be a linear symmetric operator on the Hilbert space ${X}_0$ with its domain $D(\mathcal{L}_0)$ being dense in ${X}_0$. Assume that there exists a constant $c>0$ such that
	\begin{eqnarray*}
		(\mathcal{L}_0v,v)_{{X}_0}\geq c||v||^2_{{X}_0}\;\;\mbox{for all}\;\;v\in D(\mathcal{L}_0).
	\end{eqnarray*}
	Then $\mathcal{L}_0$ has a self-adjoint extension $\mathcal{L}:D(\mathcal{L})\subseteq{X}_1\subseteq{X}_0\rightarrow{X}_0$where ${X}_1$ denotes the {energetic} Hilbert space endowed with the {energetic} scalar product $(v,w)_{{X}_1}=(\mathcal{L}v,w)_{{X}_0}$ for all $u,v\in{X}_1$ and the {energetic} norm $||v||_{{X}_1}^2=(\mathcal{L}v,v)_{{X}_0}$. Furthermore, the operator equation
	\begin{eqnarray*}
		\mathcal{L}v=f, f\in {X}_0,
	\end{eqnarray*}
	has a unique solution $v\in D(\mathcal{L})$. In addition, if $\hat{\mathcal{L}}:{X}_1\rightarrow{X}_1^*$ denotes the {energetic} extension of $\mathcal{L}$, then $\hat{\mathcal{L}}$ is the canonical isomorphism from ${X}_1$ to its dual ${X}_1^*$ and the operator equation
	\begin{eqnarray*}
		\hat{\mathcal{L}}v=f, f\in {X}_1^*,
	\end{eqnarray*}
	has also a unique solution $v\in{X}_1$.
\end{theorem}
With Theorem \ref{T1} in hand, we discuss the left-hand side of Eq. \eqref{eop3}.
It is well known that Theorem \ref{T1} is applicable to the operator $\mathcal{L}_0:D(\mathcal{L}_0)\subseteq  L^2(-L,L)\rightarrow L^2(-L,L)$,
$\mathcal{L}_0Q=-Q''(z)$,  with domain of definition, $D(\mathcal{L}_0)$, the space of $C^{\infty}$-functions on $(-L,L)$. Since $D(\mathcal{L}_0)$ is dense in $\mathcal{X}_0$, and inequality (\ref{e6}) holds, the Friedrichs extension of $\mathcal{L}_0$ is the operator $\mathcal{L}:D(\mathcal{L})\rightarrow \mathcal{X}_0$ where
\begin{eqnarray*}
 D(\mathcal{L})=\left\{v\in\mathcal{X}_1\;:\;\mathcal{L}v\in L^2(-L,L)\right\}.
\end{eqnarray*}
Consequently, the equation
\begin{eqnarray}
 \label{e7}
-Q''(z)=f,\;\;\mbox{for every}\;\;f\in L^2(-L,L),
\end{eqnarray}
has a unique solution in $D(\mathcal{L})$. Thus, we shall consider the right-hand side of Eq. \eqref{eop3} as a suitable mapping on $L^2(-L,L)$.  For its linear part, we have the following lemma proved in \cite[Proposition 1, pg.\,268]{Smets}.
\begin{lemma}
	\label{L1}
	The linear operators
	$$A_{1}[Q(z)]=Q(z+1)-Q(z)=\int_{z}^{z+1}Q'(s)ds,\;\;A_{2}[Q(z)]=Q(z)-Q(z-1)=\int_{z-1}^{z}Q'(s)ds,$$
	are continuous from $\mathcal{X}_1$ to $ L^2(-L,L)\cap L^{\infty}(-L,L)$ and $||A_{i}Q||_{L^\infty}\leq ||Q||_{{\mathcal{X}_1}}$, $||A_iQ||_{{\mathcal{X}_0}}\leq ||Q||_{{\mathcal{X}_1}}$, $i=1,2$.
\end{lemma}
Using Theorem \ref{T1}, we treat the following auxiliary linear, non-homogeneous problem
\begin{eqnarray}
\label{eop41}
-Q''(z)=\frac{1}{c^2}\left\{\kappa\left([\Psi(z)-\Psi(z-1)]-[\Psi(z+1)-\Psi(z)]\right)+V^\prime(\Psi(z))\right\}.
\end{eqnarray}
for some arbitrary fixed $\Psi\in \mathcal{X}_1$, as an equation  of the form \eqref{e7}. In particular, we have the following result:
\begin{proposition}
\label{P11}
For any $\Psi\in \mathcal{X}_1$, the equation (\ref{eop41}) has a unique solution $Q\in D(\mathcal{L})\subset\mathcal{X}_1$.
\end{proposition}
\textbf{Proof:} Equation (\ref{eop41}) can be rewritten in the form 
\begin{eqnarray}
\label{e91}
-Q''(z)=\frac{1}{c^2}\left\{\kappa{\left([A_2[\Psi(z)]-A_1[\Psi(z)]\right)}+V^\prime(\Psi(z))\right\}:=\mathcal{F}[\Psi(z)].
\end{eqnarray}
Due to the continuous embedding  $\mathcal{X}_1\subset L^{r}(-L,L)$ for any $1\leq r\leq \infty$ and condition (\ref{eq:assumptions2}),  we have that for some constant $C=C(K,\beta)>0$,
\begin{eqnarray*}
||V^\prime(\Psi)||^2_{\mathcal{X}_0}=\int_{-L}^{L}(V^\prime(\Psi(z)))^{2}dz\leq
K^2 \int_{-L}^{L}|\Psi(z)|^{2(\beta+1)}dz
\le C||\Psi||_{\mathcal{X}_1}^{2(\beta+1)}.
\end{eqnarray*}
Then, for the right-hand side of Eq.\,\eqref{e91} we get
\begin{eqnarray}
\label{eop51}
||\mathcal{F}[\Psi]||_{\mathcal{X}_0}&\leq& \frac{1}{c^2}\left\{\kappa||A_2[\Psi]-A_1[\Psi]||_{\mathcal{X}_0} +||V^\prime(\Psi(z))||_{{\mathcal{X}_0}}\right\}\nonumber\\
&\leq&\frac{1}{c^2}\left\{2\kappa||\Psi||_{\mathcal{X}_1}+C ||\Psi||_{\mathcal{X}_1}^{2(\beta+1)}\right\}.
\end{eqnarray}
Thus, 
${  \mathcal{F}[\Psi] \in L^2(-L,L)}$, and by virtue of Theorem \ref{T1} 
Eq.\,(\ref{e91})
has a unique solution  $Q\in D(\mathcal{L})$. $\square$\vspace{0.2cm}

We now proceed with the implementation of the fixed point argument. To this aim we consider for some $R>0$ the closed ball of $\mathcal{X}_1$, $\mathcal{B}_R:=\left\{\psi\in \mathcal{X}_{1}\,:\,||\psi||_{{\mathcal{X}_1}}\leq R\right\}$. Proposition \ref{P11} shows that the map $\mathcal{T}:\mathcal{X}_1\rightarrow \mathcal{X}_1$ defined as
\begin{eqnarray*}
 \mathcal{T}[\Psi]=Q,
\end{eqnarray*}
where $Q$ is the unique solution of the auxiliary problem (\ref{e91}), is well defined. Hence we introduce $\Psi_1,\Psi_2\in\mathcal{B}_R$ such that
$Q=\mathcal{T}[\Psi_1]\;\;\mbox{and}\;\;P=\mathcal{T}[\Psi_2].$
Then the difference $Y=Q-P$ satisfies the equation
\begin{eqnarray}
 \label{e141}
-Y''(z)&=&\mathcal{F}[\Psi_1(z)]-\mathcal{F}[\Psi_2(z)]\nonumber\\
&=&\frac{1}{c^2}\left\{\kappa\left(A_1[\Psi_2(z)]-A_1[\Psi_1(z)]+A_2[\Psi_1(z)]-A_2[\Psi_2(z)]\right)\right\}\nonumber\\
&+&\left.V^\prime(\Psi_1(z))-V^\prime(\Psi_2(z))\right\}.\end{eqnarray}
From Lemma \ref{L1} the linear operators $A_i:\mathcal{X}_1\rightarrow L^2(-L,L)\cap L^{\infty}(-L,L)$ are globally Lipschitz,
\begin{eqnarray}
 \label{e15}
||A_i\psi_1-A_i\psi_2||_{\mathcal{X}_0}\leq||\psi_1-\psi_2||_{\mathcal{X}_1},\\
||A_i\psi_1-A_i\psi_2||_{L^\infty}\leq||\psi_1-\psi_2||_{\mathcal{X}_1},\,\,\, {{{i=1,2.}}}\nonumber
\end{eqnarray}
To estimate the difference of the remaining terms in Eq. \eqref{e141}, we use \eqref{eq:assumptions2}	 and  the embedding $\mathcal{X}_1\subset L^{\infty}(-L,L)$ with its embedding constant $C_*$,
\begin{eqnarray}
\label{emv161}
\int_{-L}^{L}|V^\prime(\Psi_1(z))-V^\prime(\Psi_2(z))|^2dz&\leq& K^2\int_{-L}^{L}(|\Psi_1(z)|^\beta+|\Psi_2(z)|^\beta)^2|\Psi_1(z)-\Psi_2(z)|^2dz\nonumber\\
&\le& K^2\max_{-L\le z \le L}(|\Psi_1(z)|^\beta+|\Psi_2(z)|^\beta)\int_{-L}^{L}|\Psi_1(z)-\Psi_2(z)|^2dz\nonumber\\
&=&2K^2 C_*^{2\beta} R^{2\beta}||\Psi_1(z)-\Psi_2(z)||_{\mathcal{X}_0}^2.
\end{eqnarray}
Hence, for the right-hand  side of (\ref{e141}) we get
\begin{eqnarray}
 \label{e17}
||\mathcal{F}[\Psi_1(z)]-\mathcal{F}[\Psi_2(z)]||_{\mathcal{X}_0}\leq M\||\Psi_1-\Psi_2||_{\mathcal{X}_1},
\end{eqnarray}
where the constant $M$ is given by 
\begin{eqnarray}
\label{e17b1}
M=\frac{2}{c^2}\left(\kappa+KC_*R^\beta\right).
\end{eqnarray}
Next, by multiplying (\ref{e141}) in the $ L^2(-L,L)$-scalar product and using the  Cauchy-Schwarz inequality  and Young's inequality,  we get the estimate
\begin{eqnarray}
\label{e18}
 ||Y||_{\mathcal{X}_1}^2&\leq& ||\mathcal{F}[\Psi_1(z)]-\mathcal{F}[\Psi_2(z)]||_{\mathcal{X}_0}\,||Y||_{\mathcal{X}_0}\nonumber\\
&\leq&M\sqrt{C(L)}||\Psi_1-\Psi_2||_{\mathcal{X}_1}\,||Y||_{\mathcal{X}_1}\nonumber\\
&\leq&\frac{1}{2}||Y||_{\mathcal{X}_1}^2
+\frac{M^2C(L)}{2}||\Psi_1-\Psi_2||_{\mathcal{X}_1}^2.
\end{eqnarray}
Note that the Poincar\'{e} inequality (\ref{e6}) has been used. From (\ref{e18}) we derive
\begin{eqnarray}
\label{e19}
||Y||_{\mathcal{X}_1}^2=||\mathcal{T}[\Psi_1]-\mathcal{T}[\Psi_2]||_{\mathcal{X}_1}^2\leq M^2C(L)||\Psi_1-\Psi_2||_{\mathcal{X}_1}^2.
\end{eqnarray}
We conclude that if the Lipschitz constant satisfies
\begin{eqnarray}
\label{e20}
M\sqrt{C(L)}<1,
\end{eqnarray}
then the map $\mathcal{T}:\mathcal{B}_R\rightarrow\mathcal{B}_R$ is a contraction.
Hence the map $\mathcal{T}$ satisfies the assumptions of the Banach Fixed Point Theorem and has a unique fixed point. By the assumptions we have that $\mathcal{F}(0)=0$. Therefore we deduce that if (\ref{e20}) holds, then the unique fixed point is the trivial one. 
\emph{Thus, nontrivial solutions exist only if (\ref{e20}) is violated, that is, when} 
\begin{eqnarray}
\label{e20a1}
M\sqrt{C(L)}>1.
\end{eqnarray}

Regarding the upper bound for the velocity we summarise in
\setcounter{theorem}{4}
\begin{theorem}
 \label{T21}
An upper bound for the velocity $c$ of nontrivial periodic TWSs      $q_n(t)=Q(n-ct)=Q(z)$ of prescribed norm $||Q||_{\mathcal{X}_1}\leq R$ to the system (\ref{eq:system}), on the periodic lattice $-L\le n \le L$, is given by
\begin{eqnarray}
 \label{e211}
c^2<2(\kappa+KC_*R^\beta)C(L).
\end{eqnarray}
\end{theorem}

\paragraph
{\it Remarks on the physical significance of the estimates for the TWSs and their velocity.} The estimates on the TWSs proved in Theorems \ref{ThII.1}, \ref{ThII.2} and \ref{T21} implicate a coherent dependence on the lattice parameters, the frequency $\Omega$ and $R$ and the velocity $c$. Motivated by the discussion of \cite[Section A, pg. 9]{DNBRS}, we aim to discuss the potential physical relevance of these estimates:
\begin{enumerate}
\item For fixed  $\overline{m}$ and $\kappa$, we observe that 
\begin{eqnarray*}
\label{lim1}
\lim_{\Omega\rightarrow\infty}\mathscr{R}_{max}&=&\lim_{\Omega\rightarrow\infty}\mathscr{R}_{crit}=\infty,\;\;\mbox{for fixed $c$},\\
\label{lim2}
\lim_{c\rightarrow\infty}\mathscr{R}_{max}&=&\lim_{c\rightarrow\infty}\mathscr{R}_{crit}=\infty,\;\;\mbox{for fixed $\Omega$}.
\end{eqnarray*}
Both limits are physically relevant in the sense that in the limit of arbitrary large frequency or velocity,  a type of  ``energy'' of the solution, measured herein in the norm of $\mathcal{X}_0$, should become also arbitrarily large. This behavior can be relevant to energy localization phenomena \cite{His} or  the notion of quasi-collapse \cite{Kim}. Note that  in the second limit as $c\rightarrow\infty$, the growth of the norm in $\mathcal{X}_0$, implies due to the Poincar\'{e} inequality \eqref{e6} (or due to Sobolev embeddings), the growth of the kinetic energy of the TWSs, which is consistent with the growth of $c$.
\item The above coherent dependence on the lattice parameters is also evident in the derived upper bound for the velocity \eqref{e211}. Theorem \ref{T21} justifies that TWs of given ``energy'' measured by the norm of $\mathcal{X}_1$ can evolve with velocity satisfying the upper-bound \eqref{e211}. 
\end{enumerate}	

\section{Soft on-site potentials}\label{section:soft}
\setcounter{equation}{0}
In this section we study the KG lattice with soft  on-site potentials of the form
 \begin{equation}
  V(x)=-\frac{\omega_0^2}{2}x^2+\frac{a}{p+1}x^{p+1},\,\,\,a>0,\,\,p>1,\label{eq:DWP}
 \end{equation}
 where we assume that $p+1=2r>0$, with integer $r$. That is, $p$ is an odd integer so that  the on-site potentials $V(x)$ possess the reflection symmetry $V(x)=V(-x)$.
 The standard quartic double-well potential is obtained for $p=3$.
We prove the existence of periodic TWSs on finite lattices with imposed periodic boundary conditions
 utilising the Mountain Pass Theorem.

 \subsection{Periodic TWSs: Existence by the Mountain Pass Theorem}
 
We  consider periodic TWSs satisfying (\ref{eq:zeromean}), that is,  solutions $Q(z)$ performing oscillations about $Q=0$ so that the associated energy $E> 0$. Hence, the  Poincar\'{e} inequality applies.
These   periodic TWSs, as solutions of (\ref{eq:TW}), are critical points of the action functional  $S:\,\mathcal{X}_1 \rightarrow {\mathbb{R}}$ given by
  \begin{equation}
   S(Q)=
   %\int_{-L}^{L}\, L(Q(z),Q^\prime(z))du=
   \int_{-L}^{L}\left[\frac{c^2}{2}\left(Q^\prime(z)\right)^2+\frac{\omega_0^2}{2}Q^2(z)-\frac{a}{p+1}Q^{p+1}(z)-
    \frac{\kappa}{2}[Q(z+1)-Q(z)]^2\right]dz.\nonumber
   % \label{eq:actiondw}
  \end{equation}

  We get
  \begin{eqnarray}
   <S^\prime (Q),P>&=& \int_{-L}^{L}\left[c^2Q^\prime(z) P^\prime(z)+\omega_0^2 Q(z) P(z)-a Q^p(z)P(z)\right.\nonumber\\
   &+&\left.
   \kappa[Q(z+1)-2Q(z)+Q(z-1)]P(z)\right]dz,\,\,\, P,Q \in \mathcal{X}_1,\label{eq:derivative}
  \end{eqnarray}
where  $<\cdot,\cdot>$ is the standard duality bracket between  $\mathcal{X}_1$ and its dual  $\mathcal{X}_1^*$; by the definition of the derivative $S'$, for any $Q\in\mathcal{X}_1$, the functional $S'(Q):\mathcal{X}_1\rightarrow \mathbb{R}$ is a linear functional acting on any $P\in\mathcal{X}_1$ as $S'(Q)[P]= <S^\prime (Q),P>$ \cite{Chow}.   

To prove the existence of TWSs we facilitate the Mountain Pass Theorem (MPT).
  We recall [\cite{Chow}, Definition 4.1, p. 130] (Palais-Smale (PS) condition) and [\cite{Chow}, Theorem 6.1, p. 140] (Mountain Pass Theorem (MPT)
of Ambrosetti-Rabinowitz \cite{Ambrosetti}).
\begin{definition}
Let $X$ be a Banach space and ${\bf E}:X\rightarrow \mathbb{R}$ be $C^1$. We say that ${\bf E}$ satisfies condition (PS) if, for any sequence $\{u_n\} \in X$ such that $|{\bf E}(u_n)|$ is bounded and ${\bf E}^\prime (u_n) \rightarrow 0$ as $n \rightarrow \infty$ in $X^*$, there  exists a convergent subsequence. 
\end{definition}
\begin{theorem}
 Let ${\bf E}:X\rightarrow \mathbb{R}$ be $C^1$ and satisfy (A) ${\bf E}(0)=0$, 
(B) $\exists \rho >0$, $\alpha>0$: $||u||_X=\rho$ implies ${\bf E}(u)\ge \alpha$, (C) $\exists u_1 \in X$: $||u_1||_X\ge \rho$ and ${\bf E}(u_1)<\alpha$.
Define 
\begin{equation}
 \Gamma =\left\{\gamma \in C^0([0,1],X)\,:\,\gamma(0)=0,\, \gamma(1)=u_1\right\}.
\end{equation}
Let $F_{\gamma}=\{\gamma(t)\in X\,:\,0\le t\le 1\}$ and $\mathcal{L}=\{F_\gamma\,:\,\gamma\in \Gamma\}.$ If ${\bf E}$ satisfies (PS), then 
\begin{equation}
 \beta:=\inf_{F_\gamma \in \mathcal{L}}\,\sup \{{\bf E}(v)\,:\,v\in F_\gamma\}\ge \alpha,
\end{equation}
is a critical value of the functional ${\bf E}$.

 \end{theorem}

 We proceed by proving the validity of the assumptions of the MPT:
 
\begin{lemma}
\label{PS}
Assume either\\ 
(i) $\omega_0^2\ge4\kappa$,\\
or\\
(ii) $\omega_0^2<4\kappa$ and $c^2>C(L)(4\kappa-\omega_0^2)$.

Then the  functional $S$ satisfies the PS-condition.
\end{lemma}

 \noindent{\bf Proof:}  Assume $S(Q_m)$ is bounded, i.e. $|S(Q_m)|\le M$ for all $m\in {\mathbb{N}}$, and 
 \begin{equation}
 	\label{clar1}
 S^\prime(Q_m)\rightarrow 0,\;\;\mbox{as $m\rightarrow \infty$ in $\mathcal{X}_1^*$}.
\end{equation} 
Recall also that 
\begin{equation}
	\label{clar2}
|<S^\prime (Q_m),Q_m>|\leq ||S^\prime (Q_m)||_{\mathcal{X}_1^*}||Q_m||_{X_1},	
\end{equation}
by the standard inequality for the duality bracket \cite[Sec. 21.5, pg. 251]{Zeidler2}. 
Then, using \eqref{clar1},we deduce that for any $\epsilon >0$,  there exists $N(\epsilon) \in \mathbb{N}$ such that 
\begin{equation}
\label{clar3}
||S'(Q_m)||_{\mathcal{X}_1^*}\leq \epsilon,\;\;\mbox{for $m>N(\epsilon)$}.
\end{equation}	 
Applying \eqref{clar3} for $\epsilon\leq 1$ and using the inequality \eqref{clar2}, we deduce that 
\begin{equation}
	\label{clar4}
|<S^\prime (Q_m),Q_m>|\leq ||Q_m||_{\mathcal{X}_1},\;\;\mbox{for $m>N(\epsilon)$}.	
\end{equation}
Hence, for chosen $b\in (1/(p+1),1/2)$, we deduce from \eqref{clar4}, that 
\begin{equation}
\label{clar5}
b \,|<S^\prime (Q_m),Q_m>|\leq ||Q_m||_{\mathcal{X}_1}<||Q_m||_{\mathcal{X}_1}+1,\;\;\mbox{for $m>N(\epsilon)$}.	
\end{equation}
We will use \eqref{clar5} to prove that $Q_m$ is bounded in $\mathcal{X}_1$, as follows: First, we derive the inequality
 \begin{eqnarray}
  1+M+|| Q_m||_{\mathcal{X}_1}&\ge&M-b <S^\prime (Q_m),Q_m\ge 
  S(Q_m)-b <S^\prime (Q_m),Q_m>\nonumber\\
  &=&\left(\frac{1}{2}-b\right)\int_{-L}^{L}\left[c^2\left(Q^\prime_m(z)\right)^2+\omega_0^2Q_m^2(z)-
   \kappa[Q_m(z+1)-Q_m(z)]^2\right]dz\nonumber\\
   &-&\frac{a}{p+1}(1-(p+1)b)\int_{-L}^{L}Q_m^p(z)dz\nonumber\\
   &\ge&   \left(\frac{1}{2}-b\right)\left[c^2 || Q^\prime_m||_{L^2}^2 +(\omega_0^2-4\kappa) || {Q}_m||_{L^2}^2\right].
 %   - \frac{a}{4}(1-4b)|| Q_||_{L^4^4.
   \end{eqnarray}
If assumption (i) holds, then 
\begin{equation}
 1+M+|| Q_m||_{\mathcal{X}_1}\ge \left(\frac{1}{2}-b\right)c^2 || {Q}_m||_{\mathcal{X}_1}^2,
\end{equation}
and if assumption (ii) holds, then 
\begin{equation}
 1+M+|| Q_m||_{\mathcal{X}_1}\ge \left(\frac{1}{2}-b\right)\left(c^2 -C(L)(4\kappa-\omega_0^2)\right]|| {Q}_m||_{\mathcal{X}_1}^2,
\end{equation}
implying that $Q_m$ is bounded in $\mathcal{X}_1$. Hence, there is a subsequence of $Q_m$ (not relabeled) and a $Q\in \mathcal{X}_1$ such that 
  $Q_{m}\rightarrow Q$ weakly in $\mathcal{X}_1$, so that by Sobolev compact embedding, one has the strong convergence $Q_{m} \rightarrow Q$ in $L^2(-L,L)$ (and in $C([-L,L])$). Using H\"older's inequality and the embedding $\mathcal{X}_1 \subset L^\infty(-L,L)$, we obtain  
  \begin{eqnarray}
   || Q_m-Q||_{\mathcal{X}_1}^2&=&\frac{1}{c^2}\left(<S^\prime (Q_m)-S^\prime (Q),Q_m-Q>\right.\nonumber\\
   &-&
   \frac{1}{c^2}\int_{-L}^{L}\left[\left(\omega_0^2(Q_m(z)-Q(z))^2-\kappa[Q_m(z+1)-Q(z+1)-(Q_m(z)-Q(z))]^2\right.\right.\nonumber\\
 &-&\left.\left.a(Q_m(z)-Q(z))^{p+1}\right]dz\right)\nonumber\\
 &\le& \frac{1}{c^2}|<S^\prime (Q_m)-S^\prime(Q),Q_m-Q>|\nonumber\\
 &+&\frac{1}{c^2}\left((\omega_0^2+4\kappa)||Q_m-Q||_{L^2}+a ||Q_m^p-Q^p||_{L^2(-L,L)}\right) ||Q_m-Q||_{L^2(-L,L)}.\label{eq:PSconvergence1}
 \end{eqnarray}
The first term on the right-hand side of (\ref{eq:PSconvergence1}) converges to zero because by assumption
$<S^\prime (Q_m)-S^\prime(Q),Q_m-Q>\rightarrow 0$ as 
$m\rightarrow \infty$.
The last converges to zero  by strong convergence. Thus, $||Q_m-Q||_{\mathcal{X}_1}=0$ so that $(Q_m)_{m\in{\mathbb{Z}}}$ has a strongly convergent subsequence and the proof is finished.   \ $\square$

\begin{lemma}
	\label{C1}
	The functional $S$ is $C^1$ on $\mathcal{X}_1$.
	\end{lemma}
 
 \noindent{\bf Proof:} The functional $S$ can be expressed as 
 \begin{equation}
  S(q)=\frac{c^2}{2}(Q,Q)+\Gamma(Q),
 \end{equation}
 where
 \begin{equation}
  \Gamma(Q)=\int_{-L}^{L}\left[\frac{\omega_0^2}{2}Q^2(z)-\frac{a}{p+1}Q^{p+1}(z)-
  \frac{\kappa}{2}(Q(z+1)-Q(z))^2\right]dz.
 \end{equation}
 Continuity of the quadratic term $(Q,Q)$ is obvious. 
By using The embedding $\mathcal{X}_1\subset L^\infty(-L,L)$ which  implies that $Q^p\in \mathcal{X}_0$ and the Poincar\'{e} inequality \eqref{e6} we get the estimate,
 \begin{equation}
  ||Q^{p}||^2_{\mathcal{X}_0}=\int_{-L}^L |Q(z)|^{2p}dz\le C(L) ||Q||_{\mathcal{X}_1}^{2p}.
 \end{equation}
Then, we have that 
 \begin{equation}
 \label{eqX0}
  \left|\Gamma(Q)\right|\le \frac{\omega_0^2}{2}||Q||_{L^2}^2+\frac{a}{p+1}C(L) ||Q||_{\mathcal{X}_1}^{p+1}+2\kappa ||Q||_{L^2}^2\le C (L)\left( \left(\frac{\omega_0^2}{2}+2\kappa\right)||Q||_{\mathcal{X}_1}^2+\frac{a}{p+1}||Q||_{\mathcal{X}_1}^{p+1}\right)<\infty.
 \end{equation}

 The Gateaux derivative of $\Gamma$ exists and is given by
 \begin{eqnarray}
   <\Gamma^\prime (Q),h>&=& \int_{-L}^{L}\left[\omega_0^2 Q(z)-a Q^p(z)\right.\nonumber\\
   &+&\left.
   \kappa(A_1[Q(z)]-A_2[Q(z)])\right]h(z)dz.\label{eq:derivativeG}
 \end{eqnarray}

 To prove that $\Gamma^\prime$ is continuous, we let $||h||_{\mathcal{X}_1}\le 1$ and $Q_m \rightarrow Q$ in $\mathcal{X}_1$. Then 
 \begin{eqnarray}
  \left|<\Gamma^\prime(Q_m)-\Gamma^\prime(Q),h>\right|&=&
  \left|\int_{-L}^{L}\left[\omega_0^2(Q_m(z)-Q(z))-a(Q_m^{p}(z)-Q^{p}(z))\right.\right.\nonumber\\
  &+&\left.\left.
  \kappa(A_1[Q_m(z)]-A_2[Q_m(z)]-(A_1[Q(z)]-A_2[Q(z)]))\right]h(z)dz\right|\nonumber\\
  &\le& \omega_0^2||Q_m-Q||_{L^2}^2+a C\left(||Q_m||_{L^\infty},||Q||_{L^\infty}\right)||Q_m-Q||_{L^2(-L,L)}^2+2\kappa ||Q_m-Q||_{L^2(-L,L)}^2\nonumber\\
  &=& \left(\omega_0^2+aC\left(||Q_m||_{L^\infty},||Q||_{L^\infty}\right)+2\kappa\right)||Q_m-Q||_{L^2(-L,L)}^2\le C \epsilon,\nonumber
 \end{eqnarray}
 for $m$ sufficiently large. Hence
 \begin{equation}
  \left|<\Gamma^\prime(Q_m)-\Gamma^\prime(Q),h>\right|\rightarrow 0\,\,\, {\rm as}\,\,\,m\rightarrow \infty,\nonumber
 \end{equation}
and the proof of the lemma is completed.\ \ $\square$

Obviously, $S(0)=0$, therefore condition (A) of the MPT is satisfied. 
For the remaining conditions (B) and (C), the proofs are given as follows:

\textit{ Proof of (B):}\qquad We distinguish the cases (i) $\omega_0^2\ge 4\kappa$, and (ii)  $\omega_0^2<4\kappa$, $c^2> C(L)(4\kappa-\omega_0^2)$.

(i) With the aid of the estimate
\begin{equation}
   \int_{-L}^{L}\left[\frac{a}{p+1}Q^{p+1}(z)+
    \frac{\kappa}{2}[Q(z+1)-Q(z)]^2\right]dz  \le  \frac{a}{p+1}C(L)||Q||_{\mathcal{X}_1}^{p+1}    +2\kappa ||Q||_{L^2(-L,L)}^2,\nonumber                                                                                      
\end{equation}
we get the sufficient condition
\begin{eqnarray}
 \frac{c^2}{2}||Q||_{\mathcal{X}_1}^2>  \frac{a}{p+1}C(L) ||Q||_{\mathcal{X}_1}^{p+1},\nonumber                                         
\end{eqnarray}
for $S(Q)$ being positive. 
Conclusively, for $||Q||_{\mathcal{X}_1}$ small enough, say $||Q||_{\mathcal{X}_1}=\rho$, and 
\begin{equation}
0<\rho<\left(\frac{(p+1)c^2}{2a C(L)}\right)^{1/{(p-1)}} 
\end{equation}
there is a $\alpha>0$ with $S(Q)\ge \alpha$ for all $||Q||_{\mathcal{X}_1}=\rho$.

(ii) In this case we use the estimate
\begin{equation}
   \int_{-L}^{L}\left[\frac{a}{p+1}Q^{p+1}(z)+
    \frac{\kappa}{2}[Q(z+1)-Q(z)]^2\right]dz  \le  \frac{a}{p+1}C(L)||Q||_{\mathcal{X}_1}^{p+1}    +2\kappa C(L)||Q||_{\mathcal{X}_1
    }^2.\nonumber  
\end{equation}
Using the  Poincar\'{e} inequality we derive the following condition for $S>0$:
\begin{eqnarray}
 \frac{1}{2}\left(c^2-(4\kappa-\omega_0^2)C(L)\right)||Q||_{\mathcal{X}_1}^2>  \frac{a}{p+1}C(L) ||Q||_{\mathcal{X}_1}^{p+1}.\nonumber                                            
\end{eqnarray}
Hence, for 
$||Q||_{\mathcal{X}_1}$ small enough, say $||Q||_{\mathcal{X}_1}=\rho$, and 
\begin{equation}
0<\rho<\left(\frac{(p+1)}{2a C(L)}(c^2-(4\kappa-\omega_0^2))\right)^{1/{(p-1)}}, 
\end{equation}
there is a $\alpha>0$ with $S(Q)\ge \alpha$ for all $||Q||_{\mathcal{X}_1}=\rho$.
\\
\textit{Proof of (C):}\qquad For $||Q||_{\mathcal{X}_0}\neq 0$ we observe that
\begin{equation}
 S(tQ)=
   \int_{-L}^{L}\left[\frac{c^2}{2}t^2\left(Q^\prime(z)\right)^2+\frac{\omega_0^2}{2}t^2Q^2(z)-\frac{a}{p+1}t^{p+1}Q^{p+1}(z)-
    \frac{\kappa}{2}t^2[Q(z+1)-Q(z)]^2\right]dz \rightarrow -\infty, 
\end{equation}
as $t\rightarrow \infty$. To summarise, by virtue of the MPT we  state
\begin{theorem}
Let either\\ 
(i) $\omega_0^2\ge4\kappa$,\\
or\\
(ii) $\omega_0^2<4\kappa$ and $c^2>C(L)(4\kappa-\omega_0^2)$.

Then the  system (\ref{eq:system}) on the periodic lattice $-L\le n \le L$ with an  on-site potential (\ref{eq:DWP}) %possessing reflection symmetry $V(x)=V(-x)$, 
has at least one nontrivial periodic TWS.
\end{theorem}

% We conclude this section with the remark that {\it positive (negative)} periodic TWSs  associated with oscillations of $Q(z)$ about the minima $\tilde{Q}_{\pm}=\pm (\omega_0^2/a)^{1/(p-1)}$ of the potential $V(Q)$, possessing reflection symmetry $V(Q)=V(-Q)$,  can also be treated using the approach above. These TWSs are characterised by oscillations of $Q$ without changes of sign and possess energy $E<0$. One only needs to consider the system (\ref{eq:TW}) expressed 
% in the shifted variable $Q\rightarrow Q\mp \tilde{Q}_{\pm}$. Note that the linear operators $A_{1,2}$ are invariant with respect to the shift operator $Q\rightarrow Q+Q_0$. 

\subsection{Thresholds for the average kinetic energy of TWSs of prescribed speed}

For the derivation of lower bounds of the average kinetic energy of TWSs 
we utilise the fixed point method outlined in Section  \ref{section:hard}.
In the present case, using Theorem \ref{T1}, we treat the following auxiliary linear, non-homogeneous problem
\begin{eqnarray}
\label{eop4}
-Q''(z)=\frac{\kappa}{c^2}\left\{\left(\Psi(z)-\Psi(z-1)\right)-\left(\Psi(z+1)-\Psi(z)\right)-\omega_0^2\Psi(z)+a\left(\Psi(z)\right)^p\right\},
\end{eqnarray}
for some arbitrary fixed $\Psi\in \mathcal{X}_1$, as an equation  of the form \eqref{e7}. We have the following result.
\begin{proposition}
\label{P1}
For any $\Psi\in \mathcal{X}_1$, the equation (\ref{eop4}) has a unique solution $Q\in D(\mathcal{L})\subset\mathcal{X}_1$.
\end{proposition}
\textbf{Proof:} We reformulate Eq.\,(\ref{eop4}) as  
\begin{eqnarray}
\label{e9}
-Q''(z)=\frac{1}{c^2}\left\{\kappa{\left(A_2\left[\Psi(z)\right]-A_1\left[\Psi(z)\right]\right)}-\omega_0^2\Psi(z)+a\left(\Psi(z)\right)^p\right\}:=\mathcal{F}\left[\Psi(z)\right].
\end{eqnarray}
We use again the estimate (\ref{eqX0}) for $\Psi$,
\begin{eqnarray*}
||\Psi^p||^2_{\mathcal{X}_0}=\int_{-L}^{L}|\Psi(z)|^{2p}dz\leq C(L)||\Psi||_{\mathcal{X}_1}^{2p}.
\end{eqnarray*}
Then, for the right-hand side of Eq. \eqref{e9} we get
\begin{eqnarray}
\label{eop5}
||\mathcal{F}[\Psi]||_{\mathcal{X}_0}&\leq& \frac{1}{c^2}\left\{\kappa||A_2[\Psi]-A_1[\Psi]||_{\mathcal{X}_0} +C(L)\omega_0^2||\Psi||_{\mathcal{X}_1}+C_{{{0}}} a||\Psi||^{2p}_{\mathcal{X}_1}\right\}\nonumber\\
&\leq&\frac{1}{c^2}\left\{2\kappa||\Psi||_{\mathcal{X}_1}+C_1\omega_0^2||\Psi||_{\mathcal{X}_1}+C_{{{0}}}
a ||\Psi||^{2p}_{\mathcal{X}_1}\right\}.
\end{eqnarray} 
Thus, 
$\mathcal{F}[\Psi] \in L^2(-L,L)$, and due to  Theorem \ref{T1}, Eq. (\ref{e9})
has a unique solution  $Q\in D(\mathcal{L})$. \ \ $\square$\vspace{0.2cm}

We proceed along the lines in Section \ref{SECIIC} by adapting the steps for the use of the fixed point argument to the present case, that is replacing the hard potential $V(Q)$ in the expressions (\ref{e141})  by the soft  potential (\ref{eq:DWP}).

In particular, we need the following estimate 
of the difference of the power nonlinearity terms 
\begin{eqnarray}
\label{emv16}
\int_{-L}^{L}|s_1^{p}-
s_2^{p}|^2&\leq& p^2\int_{-L}^{L}
\left\{\int_0^1|\xi|^{p-1}|s_1-s_2|d\theta\right\}^2,
\end{eqnarray}
where $s_1,\;s_2\in \mathbb{R}$, and  $\xi=\theta s_1
+(1-\theta)s_2$, $\theta\in (0,1)$. Then, we have that
\begin{eqnarray}
\label{e12}
\int_{-L}^{L}\left|[\Psi_1(z)]^{p}-[\Psi_2(z)]^{p}\right|^2dz&\leq& p^2||\xi||_{L^\infty}^{2(p-1)}\int_{-L}^{L}|\Psi_1(z)-\Psi_2(z)|^2dz\nonumber\\
&\leq&p^2||\xi||_{L^\infty}^{2(p-1)}C\int_{-L}^{L}|\Psi_1'(z)-\Psi_2'(z)|^2dz,
\end{eqnarray}
for $\xi(z)=\theta\Psi_1(z)+(1-\theta)\Psi_2(z)$. For the norm $||\xi||_{L^\infty}$ we have the estimate
\begin{eqnarray}
\label{eo16a}
 ||\xi||_{L^\infty}\leq\theta||\Psi_1||_{L^\infty}+(1-\theta)||\Psi_2||_{L^\infty}
\leq\theta C_*||\Psi_1||_{\mathcal{X}_1}+(1-\theta)C_*||\Psi_2||_{\mathcal{X}_1},
\end{eqnarray}
where the constant $C_*>0$ denotes the optimal constant of the embedding $\mathcal{X}_1\subset L^{\infty}(-L,L)$. Therefore, since $\Psi_1,\Psi_2\in\mathcal{B}_R$, we have that $||\xi||_{L^\infty}\leq C_*R$. Thus by using (\ref{e12}) and (\ref{eo16a}) we deduce the inequality
\begin{eqnarray*}
 ||[\Psi_1(z)]^{p}-[\Psi_2(z)]^{p}||_{\mathcal{X}_0}\leq p\sqrt{C}(C_*R)^{p-1}||\Psi_1-\Psi_2||_{\mathcal{X}_1}.
\end{eqnarray*}
Then the Lipschitz constant $M$ in the equation equivalent to (\ref{e17b1}) in \ref{subsection:upperbound} is determined by 
\begin{eqnarray}
\label{e17b}
M=\frac{1}{c^2}\left({{{2}}}\kappa+{{{C(L)}}}\omega_0^2+p\sqrt{C(L)} a C_*^{p-1}R^{p-1}\right).
\end{eqnarray}
Using again the fixed point argument we end up with the following statement:

If
 \begin{equation}
  M\sqrt{C(L)}<1
 \end{equation}
 holds, the unique fixed point is the trivial one. 
\emph{Nontrivial solutions exist only if }
\begin{eqnarray}
\label{e20a}
M\sqrt{C(L)}>1.
\end{eqnarray}
From Eq. \eqref{e17b} we derive the following condition for the  existence of nontrivial solutions:
\begin{eqnarray}
\label{fineq1}
R^2>\left[\frac{c^2-\sqrt{C(L)}({{{2}}}\kappa+{{{C(L)}}}\omega_0^2)}{C(L)}\right]^{\frac{2}{p-1}}
\left(\frac{1}{p a C_*^{p-1}}\right)^{\frac{2}{p-1}}:=T_{\mathrm{thresh}}.
\end{eqnarray}

We conclude:
\setcounter{theorem}{4}
\begin{theorem}
 \label{T2}
Consider the system (\ref{eq:system})  on the periodic lattice of $2L$ particles, $-L\leq n\leq L$ with a soft on-site potential (\ref{eq:DWP}) with reflection symmetry $V(x)=V(-x)$. Every  nontrivial periodic travelling wave solution $q_n(t)=Q(n-ct)=Q(z)$ with speed $c$ satisfying
\begin{eqnarray}
 \label{e21}
c^2>\sqrt{C(L)}({{{2}}}\kappa+{{{C(L)}}}\omega_0^2)=c_{\mathrm{crit}}^2,
\end{eqnarray}
must have average kinetic energy $T(Q)=\frac{1}{2}\int_{-L}^{L}U'(z)^2dz$ satisfying
the lower bound
\begin{eqnarray}
\label{e22}
T_{\mathrm{thresh}}<2T(Q).
\end{eqnarray}
\end{theorem}

The relation (\ref{e22}) can be regarded as a threshold value criterion for the average kinetic energy in order that  travelling waves of speed $c>c^*$ exist.

We remark that further useful quantifications of the norm and energy thresholds can be provided as explicit values of the optimal constants $C_*$ and $C_{1,*}$ of the Sobolev embeddings used in our proofs (see in \cite{bsc}).

%%%%%%%%%%%%%%%%%%%%%%%%%%%%%%%%NEW %%%%%%%%%%%%%%%%%

\vspace{0.5cm}

\noindent\textbf{Acknowledgment}
The authors are  grateful to the referees for their valuable comments and suggestions which improved considerably the presentation of the manuscript.

\vspace{0.5cm}

\noindent\textbf{Data Availability Statement}\\
The article has no associated data.\\
\\
\textbf{Authors Declarations}\\
The authors have no conflicts to disclose.\\
\\
\textbf{Authors Contributions Statement}\\
All authors contributed equally to the study conception, design and writing of the manuscript. Material preparation, data collection and analysis were performed equally by all authors.  All authors read and approved the final manuscript.
%\section*{Acknowledgment}
%We would like to thank the reviewer(s) for valuable comments and suggestions. 
%%

\end{document}